# Smarandache idempotents in certain types of group rings

**Parween Ali Hummadi**

College of Science Education

University of Salahaddin

Erbil – Kurdistan Region- Iraq

**Shadan Abdulkadr Osman**

College of Science

University of Salahaddin

Erbil- Kurdistan Region - Iraq

**Abstract:**

In this paper we study S-idempotents of the group ring $\mathbb{Z}_2 G$ where $G$ is a finite cyclic group of order $n$. We give a condition on $n$ such that every nonzero idempotent element of the group ring $\mathbb{Z}_2 G$ is Smarandache idempotent and we find Smarandache idempotents of the group ring $\mathcal{K}G$, where $\mathcal{K}$ is an algebraically closed field of characteristic $0$ and $G$ is a finite cyclic group.

**Keywords:** Idempotent, S-idempotent, group ring, algebraically closed field.

## Introduction:

Smarandache idempotent element in rings introduced by Vasantha Kandasamy [1]. A Smarandache idempotent (S-idempotent) of the ring $\mathcal{R}$ is an element $0 \neq x \in \mathcal{R}$ such that
1) $x^2 = x$
2) There exists $a \in \mathcal{R} \setminus \{0, 1, x\}$
   i) $a^2 = x$ and
   ii) $xa = a$ $(ax = a)$ or $ax = x$ $(xa = x)$.

She introduced many Smarandache concepts [2]. Vasantha Kandasamy and Moon K. Chetry discuss S-idempotents in some type of group rings [3],. A prime number $p$ of the form $p = 2^k - 1$ where $k$ is a prime number called Mersenne prime [4]. In section one of this paper we study S-idempotents of the group ring $\mathbb{Z}_2 G$ where $G$ is a finite cyclic group of order $n$. If $n = 2p$, $p$ is a Mersenne prime, we show that every nonzero idempotent element is S-idempotent and we find the number of S-idempotent element. In section two we study S-idempotents of the group ring $\mathcal{K}G$ where $\mathcal{K}$ is an algebraically closed field of characteristic $0$ and $G$ is a finite cyclic group, we show that every non trivial idempotent is S-idempotent.

## 1. S-idempotents of $\mathbb{Z}_2 G$

In this section we study S-idempotents in the group ring $\mathbb{Z}_2 G$ where $G$ is a finite cyclic group of order $n$, specially where n=2p, p is a Mersenne prime (i.e. $p = 2^k - 1$ for some prime $k$).

**Theorem 1.1.**

The group ring $\mathbb{Z}_2 G$ where G=$\langle g \mid g^m = 1 \rangle$ is a cyclic group of an odd order $m > 1$, has at least two non trivial idempotent elements, moreover no non trivial idempotent element is S-idempotent.

**Proof:** Consider the element

$$\alpha = g + g^2 + g^3 + \ldots + g^{\frac{m-1}{2}} + g^{\frac{m-1}{2}+1} + \ldots$$

$+ g^{m-1}$, of $\mathbb{Z}_2 G$. Since the coefficient of each $g^i$, $i = 1, \ldots, m$ is in $\mathbb{Z}_2$, $\alpha^2 = g^2 + g^4 + \ldots + g^{m-1} + g + g^3 + \ldots + g^{m-2}$. Hence $\alpha^2 = \alpha$, that is $\alpha$ is an idempotent element, so $(1 + \alpha)$ is also an idempotent element. It remains to show that no idempotent element of $\mathbb{Z}_2 G$ is an S-idempotent. Suppose

$$\alpha = a_1 + a_2 g + a_3 g^2 + \ldots + a_{\frac{m-1}{2}+1} g^{\frac{m-1}{2}} +$$

$\ldots + a_m g^{m-1}$, is a non trivial S-idempotent. Thus $\alpha$ is different from 0 and 1, moreover there exists $\beta$ in $\mathbb{Z}_2 G \setminus \{0, 1, \alpha\}$ such that $\beta^2 = \alpha$, let $\beta = b_1 + b_2 g + b_3 g^2 + \ldots +$

$$b_{\frac{m-1}{2}+1} g^{\frac{m-1}{2}} + \ldots + b_m g^{m-1},$$

where $b_i \in \mathbb{Z}_2$. But $\alpha^2 = \alpha$, which means that

$$a_1 + a_2 g^2 + a_3 g^4 + \ldots + a_{\frac{m-1}{2}+1} g^{m-1} +$$

$+ \ldots + a_m g^{m-2} = b_1 + b_2 g^2 + b_3 g^4 + \ldots + b_{\frac{m-1}{2}+1} g^{m-1} + \ldots + b_m g^{m-2}$.

It follows that $a_i = b_i$ for each $(1 \leq i \leq m)$. Therefore $\alpha = \beta$, which is an obvious contradiction.

The group ring $\mathbb{Z}_2 G$, where $G$ is acyclic group of an odd order may contains more than two idempotent elements as it is shown by the following example.

**Example 1.1.**

Consider the group ring $\mathbb{Z}_2 G$ where $G = \langle g \mid g^7 = 1 \rangle$ is a cyclic group of order 7. By Theorem 1.1, $g + g^2 + g^3 + g^4 + g^5 + g^6$ and $1 + g + g^2 + g^3 + g^4 + g^5 + g^6$ are idempotent elements, In addition $(g + g^2 + g^4)^2 = g^2 + g^4 + g$ and $(1 + g + g^2 + g^4)^2 = 1 + g^2 + g^4 + g$, so $1 + g + g^2 + g^4$ and $g + g^2 + g^4$ are idempotent elements. Therefore $\mathbb{Z}_2 G$ has more than two idempotent elements.

The proof of the following result is not difficult.

**Theorem 1.2.**

If $\alpha$ is an S-idempotent of the group ring $\mathbb{Z}_2 G$ where $G$ is a cyclic group of order $n$, then $(1 + \alpha)$ is an S-idempotent of $\mathbb{Z}_2 G$.

**Theorem 1.3.**

The group ring $\mathbb{Z}_2 G$, where $G = \langle g \mid g^{2n} = 1 \rangle$ is a cyclic group of order $2n$, $n$ is an odd prime, has at least two S-idempotents.

**Proof:** Let $\alpha = g^2 + g^4 + \cdots + g^{n-1} + g^{n+1} + \cdots + g^{2n-2}$. Thus

$\alpha^2 = g^4 + g^8 + \cdots + g^{2n-2} + g^2 + g^6 + \cdots + g^{2n-4} = \alpha$. Hence $\alpha$ is an idempotent element, so $(1 + \alpha)$ is also an idempotent element. We will show that $\alpha$ is S-idempotent, so let

$$\beta = g + g^{n+2} + g^3 + g^{n+4} + \cdots +$$
$$g^{\frac{n-1}{2}} + g^{\frac{3n+1}{2}} + \cdots + g^{n-2} + g^{2n-1}.$$

It is clear that $\beta^2 = \alpha$. We claim that $\alpha \beta = \beta$. For this purpose we describe the multiplication $\alpha \beta$ by the following array say $\mathcal{A}$:

$$\mathcal{A} = \begin{bmatrix} \boxed{g^3} & g^5 & \cdots & g^{n-2} & g^n & g^{n+2} & \cdots & g^{2n-3} & g^{2n-1} \\ \boxed{g^{n+4}} & g^{n+6} & \cdots & g^{2n-1} & g^{2n+1} & g^{2n+3} & \cdots & g^{n-2} & g^n \\ g^5 & g^7 & \cdots & g^n & g^{n+2} & g^{n+4} & \cdots & g^{2n-1} & \boxed{g} \\ g^{n+6} & g^{n+8} & \cdots & g & g^3 & g^5 & \cdots & g^n & \boxed{g^{n+2}} \\ \vdots & \vdots & \ddots & \vdots & \vdots & \vdots & \ddots & \vdots & \vdots \\ \boxed{g^{\frac{3n+1}{2}}} & g^{\frac{3n+5}{2}} & \cdots & g^{\frac{5n-9}{2}} & g^{\frac{5n-5}{2}} & g^{\frac{5n-1}{2}} & \cdots & g^{\frac{7n-11}{2}} & g^{\frac{7n-7}{2}} \\ \boxed{g^{\frac{n+3}{2}}} & g^{\frac{n+7}{2}} & \cdots & g^{\frac{3n-7}{2}} & g^{\frac{3n-3}{2}} & g^{\frac{3n+1}{2}} & \cdots & g^{\frac{5n-9}{2}} & g^{\frac{5n-5}{2}} \\ g^{\frac{3n+5}{2}} & g^{\frac{3n+9}{2}} & \cdots & g^{\frac{5n-5}{2}} & g^{\frac{5n-1}{2}} & g^{\frac{5n+3}{2}} & \cdots & g^{\frac{7n-7}{2}} & \boxed{g^{\frac{7n-3}{2}}} \\ g^{\frac{n+7}{2}} & g^{\frac{n+11}{2}} & \cdots & g^{\frac{3n-3}{2}} & g^{\frac{3n+1}{2}} & g^{\frac{5n+3}{2}+1} & \cdots & g^{\frac{5n-5}{2}} & \boxed{g^{\frac{5n-1}{2}}} \\ \vdots & \vdots & \ddots & \vdots & \vdots & \vdots & \ddots & \vdots & \vdots \\ \boxed{g^{n-2}} & g^n & \cdots & g^{2n-7} & g^{2n-5} & g^{2n-3} & \cdots & g^{n-8} & g^{n-6} \\ \boxed{g^{2n-1}} & g & \cdots & g^{n-6} & g^{n-4} & g^{n-2} & \cdots & g^{2n-7} & g^{2n-5} \\ g^n & g^{n+2} & \cdots & g^{2n-5} & g^{2n-3} & g^{2n-1} & \cdots & g^{n-6} & \boxed{g^{n-4}} \\ g & g^3 & \cdots & g^{n-4} & g^{n-2} & g^n & \cdots & g^{2n-5} & \boxed{g^{2n-3}} \\ g & g^3 & \cdots & g^{p-4} & g^{p-2} & g^p & \cdots & g^{2p-5} & \boxed{g^{2p-3}} \end{bmatrix}$$

That is $\mathcal{A} = [a_{ij}]_{(n-1)\times(n-1)}$, where $a_{ij}$ is the summand of $\alpha\beta$ which is equal to the product of the $i$th summand of $\beta$ with the $j$th summand of $\alpha$. This means $\alpha\beta = \sum_{i=1}^{n-1}\sum_{j=1}^{n-1} a_{ij}$. If we take the first and the third rows of this array we will see that $g^i$ occurs twice for each $i$ except $(i = 1, 3)$. By adding the terms of this two rows it remains only $g + g^3$ (observing that the coefficient of each $g^i$, i=1, 2, …,m is in $\mathbb{Z}_2$ ). Again by adding the second and the fourth rows in this array, according to the same argument it remains only $g^{p+2}+g^{p+4}$. Proceeding in this manner we will get the $(p-3)$th and the $(p-1)$th rows, and adding their terms it remains only $g^{2p-3} + g^{2p-1}$. Thus we get

$\alpha\beta = g + g^{n+2} + g^3 + g^{n+4} + \cdots + g^{\frac{n-1}{2}} + g^{\frac{3n+1}{2}} + \cdots + g^{n-2}+g^{2n-1} = \beta$.

Hence $\alpha$ is S-idempotent. By Theorem 1.2, $(1+\alpha)$ is also S-idempotent. This complete the proof.

**Lemma 1.4.**

In $\mathbb{Z}_2 G$, where $G = \langle g \mid g^{2p} = 1 \rangle$, $p$ is a Mersenne prime (i.e. $p = 2^k - 1$ for some prime $k$) $g^{2l} = g^{2^{k+1}l}$ and the elements of $S = \{g^{2l}, g^{2^2 l}, g^{2^3 l}, \ldots, g^{2^{k-1}l}, g^{2^k l}\}$ are distinct for each odd number $l$ less than $p$.

**Proof:** Since $2^{k+1}l - 2l = 2l(2^k - 1) = 2lp$, $2^{k+1}l \equiv 2l \pmod{2p}$, which implies that $g^{2l} = g^{2^{k+1}l}$. Now suppose that $g^{2l} = g^{2^t l}$ (for some $1 < t \leq k$). This means $2^t l \equiv 2l \pmod{2p}$, hence $(2^k - 1) \mid l(2^{t-1} - 1)$ yields either $(2^k - 1) \mid l$ or $(2^k - 1) \mid (2^{t-1} - 1)$. But $(2^k - 1) \mid l$ contradicts the hypothesis that $l < p$, and if $(2^k - 1) \mid (2^{t-1} - 1)$, hence $k < t - 1$, contradiction with $1 < t \leq k$.

**Lemma 1.5.**

If $p = 2^k - 1$ is a Mersenne prime, then $k \mid (2^k - 2)$.

**Proof:** Since $k$ is prime, according to Fermat's Little Theorem, $k \mid (2^k - 2)$.

Combining the last two lemmas we deduce that in the group ring $\mathbb{Z}_2 G$, where $G$ is a cyclic group generated by g of order $2p$, $p$ is a Mersenne prime (i.e. $p = 2^k - 1$ for some prime $k$), if $m = \frac{2^k - 2}{k}$, then $\alpha = g^2 + g^4 + \cdots + g^{p-1} + g^{p+1} + \cdots + g^{2p-2}$, can be partitioned to sum of $m$ elements say $\alpha_1, \alpha_2, \ldots, \alpha_m$ each $\alpha_i$ ($1 \leq i \leq m$) is of the form

$\alpha_i = g^{2l} + g^{2^2 l} + \ldots + g^{2^{k-1}l} + g^{2^k l}$,

where $l$ is an odd number.

**Theorem 1.6.**

Let $\mathbb{Z}_2 G$ be a group ring, where $G = \langle g \mid g^{2p} = 1 \rangle$ is a cyclic group of order $2p$, $p$ is a Mersenne prime. Then every element of the form $\alpha = g^{2l} + g^{2^2 l} + \cdots + g^{2^k l}$, is an S-idempotent ($l$ is an odd number).

**Proof:** Let $\alpha = g^{2l} + g^{2^2 l} + \cdots + g^{2^k l}$. By Lemma 1.4, all elements in $S = \{g^{2l}, g^{2^2 l}, \ldots, g^{2^k l}\}$ are distinct, moreover $g^{2l} = g^{2^{k+1}l}$. Hence $\alpha^2 = \alpha$. Now, let $\beta = g^l + g^{t_2} + g^{t_3} + \cdots + g^{t_k}$ and $x_i$, $i \geq 2$ be the smallest positive integer such that $x_i < 2p$. Thus $x_i \equiv 2^i l \pmod{2p}$, this means $x_i = 2^i l - 2pr$, for some $r \in \mathbb{Z}^+$. Define $t_i$ by

$$t_i = \begin{cases} \frac{1}{2} x_i & \text{if } \frac{1}{2} x_i \text{ is odd } (2 \leq i \leq k) \\ \frac{1}{2} x_i + p & \text{if } \frac{1}{2} x_i \text{ is even } (2 \leq i \leq k). \end{cases}$$

If $\frac{1}{2} x_i$ is odd, then $(g^{t_i})^2 = \left(g^{2^{i-1}l - pr}\right)^2 = g^{2^i l}$. Hence $\beta^2 = \alpha$. If $\frac{1}{2} x_i$ is even, then $(g^{t_i})^2 = g^{2^i l}$, and $\beta^2 = \alpha$ for each ($2 \leq i \leq k$). We will show that $\alpha \beta = \beta$. For this purpose as before we describe the multiplication $\alpha \beta$ in the following array say $\mathcal{A}$:

$$\mathcal{A} = \begin{bmatrix} g^{3l} & g^{5l} & g^{9l} & \cdots & g^{l(2^{k-2}+1)} & g^{l(2^{k-1}+1)} & \boxed{g^{l(2^k+1)}} \\ \boxed{g^{t_2+2l}} & g^{t_2+4l} & g^{t_2+8l} & \cdots & g^{t_2+2^{k-2}l} & g^{t_2+2^{k-1}l} & g^{t_2+2^k l} \\ g^{t_3+2l} & \boxed{g^{t_3+4l}} & g^{t_3+8l} & \cdots & g^{t_3+2^{k-2}l} & g^{t_3+2^{k-1}l} & g^{t_3+2^k l} \\ g^{t_4+2l} & g^{t_4+4l} & \boxed{g^{t_4+8l}} & \cdots & g^{t_4+2^{k-2}l} & g^{t_4+2^{k-1}l} & g^{t_4+2^k l} \\ \vdots & \vdots & \vdots & \ddots & \vdots & \vdots & \vdots \\ g^{t_{k-1}+2l} & g^{t_{k-1}+4l} & g^{t_{k-1}+8l} & \cdots & \boxed{g^{t_{k-1}+2^{k-2}l}} & g^{t_{k-1}+2^{k-1}l} & g^{t_{k-1}+2^k l} \\ g^{t_k+2l} & g^{t_k+4l} & g^{t_k+8l} & \cdots & g^{t_k+2^{k-2}l} & \boxed{g^{t_k+2^{k-1}l}} & g^{t_k+2^k l} \end{bmatrix} = [a_{ij}]_{k\times k},$$

where $a_{ij}$ is the summand of $\alpha\beta$ which is equal to the product of the $i$th summand of $\beta$ with $j$th summand of $\alpha$. This means $\alpha\beta = \sum_{i=1}^{k}\sum_{j=1}^{k} a_{ij}$. We complete the proof by the following three steps.

**Step 1:** Considering the first and the $k$th column in this array we claim that
$$a_{1j} = a_{(j+1)k} \quad \ldots(1),$$
for each $(1 \leq j \leq k-1)$, equivalently $g^{(2^j+1)l} = g^{t_{j+1}+2^k l}$.

Let $\omega = t_{j+1} + 2^k l - (2^j+1)l$. Now, $x_{j+1} \equiv 2^{j+1}l \pmod{2p}$, thus $x_{j+1} = 2^{j+1}l - 2pr$, for some $r \in \mathbb{Z}^+$. If $\frac{1}{2}x_{j+1}$ is odd, then $\frac{1}{2}x_{j+1} = 2^j l - pr$ is odd (this hold only if $r$ is odd), hence $t_{j+1} = 2^j l - pr$. So, $\omega = 2^j l - pr + 2^k l - 2^j l - l \equiv 0 \pmod{2p}$. Therefore $(2^j+1)l \equiv t_{j+1} + 2^k l \pmod{2p}$. This yields (1). If $\frac{1}{2}x_{j+1}$ is even, then $\frac{1}{2}x_{j+1} = 2^j l - pr$ is even (this hold only if $r$ is even), hence $t_{j+1} = 2^j l - pr + p$. So, $\omega = (1-r)p + lp \equiv 0 \pmod{2p}$. Hence $(2^j+1)l \equiv t_{j+1} + 2^k l \pmod{2p}$. This also yields (1). This implies that $a_{1j} + a_{(j+1)k} = 0 \pmod{2p}$, therefore by adding the terms of the first row and the $k$th column it remains only $a_{1k} = g^{l(2^k+1)}$.

**Step 2:** Consider the subarray $\mathcal{B} = (b_{ij})_{k-1 \times k-1}$ of $\mathcal{A} = (a_{ij})_{k \times k}$, where $b_{ij} = a_{(i+1)j}$ for each $(1 \leq i, j \leq k-1)$, by neglecting the first row and the $k$th column, we will show that
$$b_{ij} = b_{ji} \quad \ldots(2),$$
for all $(1 \leq i, j \leq k-1)$ such that $(i \neq j)$, equivalently $g^{t_{(i+1)}+2^j l} = g^{t_{(j+1)}+2^i l}$. Let $\omega = t_{i+1} + 2^j l - t_{j+1} - 2^i l$. Now, $x_{i+1} = 2^{i+1}l - 2pr$ and $x_{j+1} = 2^{j+1}l - 2ps$, for some $r, s \in \mathbb{Z}^+$. Thus $\frac{1}{2}x_{i+1} = 2^i l - pr$ and $\frac{1}{2}x_{j+1} = 2^j l - ps$. If $\frac{1}{2}x_{i+1}$ and $\frac{1}{2}x_{j+1}$ are even, hence $2^i l - pr$ and $2^j l - ps$ are even (this hold only if $r$ and $s$ are even), it follows $t_{i+1} = 2^i l - pr + p$ and $t_{j+1} = 2^j l - ps + p$. So, $\omega = (s-r)p \equiv 0 \pmod{2p}$. Hence $t_{i+1} + 2^j l \equiv t_{j+1} + 2^i l \pmod{2p}$. This yields (2). If $\frac{1}{2}x_{i+1}$ and $\frac{1}{2}x_{j+1}$ are odd, it is clearly $\omega = (s-r)p \equiv 0 \pmod{2p}$. Hence $t_{i+1} + 2^j l \equiv t_{j+1} + 2^i l \pmod{2p}$.

This also establishes (2). If $\frac{1}{2} x_{i+1}$ is odd and $\frac{1}{2} x_{j+1}$ is even, it is also clear that $\omega = (s - r - 1)p \equiv 0 \pmod{2p}$. Thus $t_{i+1} + 2^j l \equiv t_{j+1} + 2^i l \pmod{2p}$. This also yields (2). If $\frac{1}{2} x_{i+1}$ is even and $\frac{1}{2} x_{j+1}$ is odd, thus by using similar argument we get $t_{i+1} + 2^j l \equiv t_{j+1} + 2^i l \pmod{2p}$. This also yields (2). For all cases we get $b_{ij} + b_{ji} = 0$ $(1 \leq i, j \leq k-1)$.

**Step 3:** From Step 1 and Step 2 we get that $\alpha\beta = a_{1k} + \sum_{i=1}^{k-1} b_{ii}$ and it is not difficult to show that $\alpha\beta = \beta$ which means that $\alpha$ is an S-idempotent.

We call an S-idempotent of $\mathbb{Z}_2 G$ of the form $\alpha = g^{2l} + g^{2^2 l} + \cdots + g^{2^k l}$, where $l$ is an odd number a basic S-idempotent.

**Example 1.2.**

Consider the group ring $\mathbb{Z}_2 G$ where $G = \langle g \mid g^{62} = 1 \rangle$ is a cyclic group of order 62 (i.e. $p = 31$ and $k = 5$). By Theorem 1.7, if $l = 1$, then $\alpha = g^2 + g^4 + g^8 + g^{16} + g^{32}$ and $\beta = g + g^{33} + g^{35} + g^{39} + g^{47}$. It is clear that $\beta^2 = \alpha$. Let us describe the multiplication $\alpha\beta$ by the following array say $\mathcal{A}$:

$$\mathcal{A} = \begin{bmatrix} g^3 & g^5 & g^9 & g^{17} & g^{33} \\ g^{35} & g^{37} & g^{41} & g^{49} & g^3 \\ g^{37} & g^{39} & g^{43} & g^{51} & g^5 \\ g^{41} & g^{43} & g^{47} & g^{55} & g^9 \\ g^{49} & g^{51} & g^{55} & g & g^{17} \end{bmatrix}.$$

Hence applying Theorem 1.6, we get $\alpha\beta = g + g^{33} + g^{35} + g^{39} + g^{47} = \beta$.

**Theorem 1.7.**

If $\alpha_1$ and $\alpha_2$ are two basic S-idempotents in $\mathbb{Z}_2 G$, where $G$ is a cyclic group of order $2p$, $p$ a Mersenne prime, then $\alpha_1 + \alpha_2$ is S-idempotent.

**Proof:** Let $\alpha_1, \alpha_2$ be two distinct basic S-idempotents in $\mathbb{Z}_2 G$, so there exist $\beta_1$ and $\beta_2$ such that $\beta_1^2 = \alpha_1$, $\alpha_1 \beta_1 = \beta_1$, $\beta_2^2 = \alpha_2$ and $\alpha_2 \beta_2 = \beta_2$.

Now, $(\beta_1 + \beta_2)^2 = \beta_1^2 + \beta_2^2 = \alpha_1 + \alpha_2$, and $(\alpha_1 + \alpha_2)(\beta_1 + \beta_2) = \alpha_1 \beta_1 + \alpha_1 \beta_2 + \alpha_2 \beta_1 + \alpha_2 \beta_2 = \beta_1 + \beta_2 + \alpha_1 \beta_2 + \alpha_2 \beta_1$. We show that $\alpha_1 \beta_2 + \alpha_2 \beta_1 = 0$. By describing the multiplications $\alpha_1 \beta_2$ and $\alpha_2 \beta_1$ by the two arrays $\mathcal{A}$ and $\mathcal{B}$ respectively and using similar argument of Theorem 1.6, we get $\mathcal{A} + \mathcal{B} = 0$ that is $\alpha_1 \beta_2 + \alpha_2 \beta_1 = 0$. Therefore $\alpha_1 + \alpha_2$ is an S-idempotent.

**Theorem 1.8.**

If $\alpha_1, \alpha_2, \ldots, \alpha_n$ are $n$ basic S-idempotents in $\mathbb{Z}_2 G$ where $G$ is a cyclic group of order $2p$, $p$ is a Mersenne prime, then $\alpha_1 + \alpha_2 + \cdots + \alpha_n$ is S-idempotent.

**Proof:** Follows from Theorem 1.7.

By combining all previous results concerning the group ring $\mathbb{Z}_2 G$, where $G$ is a cyclic group of order $2p$, $p$ is a Mersenne prime we get the following result

**Theorem 1.9.**

Consider the group ring $\mathbb{Z}_2 G$ where $G$ is a cyclic group of order $2p$, $p$ is a Mersenne prime. Then
1) Every non trivial idempotent is S-idempotent.
2) The number of non trivial S-idempotents is $2(2^m - 1)$, where $m = \frac{p-1}{k}$.

**Proof:** 1) Follows from Theorems 1.6, 1.7, 1.8 and Theorem 1.2.
2) From Theorems 1.6, 1.7, and 1.8, by using the concepts of probability theory we conclude that the number of S-idempotent in $\mathbb{Z}_2 G$ is

$$\lambda = 2\left(\binom{m}{1} + \binom{m}{2} + \cdots + \binom{m}{m}\right) = 2(2^m - 1), \text{ where } m = \frac{p-1}{k}.$$

## 2. S-idempotents in the group ring of a finite cyclic group over a field of characteristic zero

In this section, we study the group ring $\mathcal{K}G$ where $\mathcal{K}$ is an algebraically closed field of characteristic 0 and $G$ is a finite cyclic group of order $n$. We get that every nontrivial idempotent element in this group ring $\mathcal{K}G$ is an S-idempotent element.

**Theorem 2.1.**

Let $\mathcal{K}$ be algebraically closed field of characteristic 0 and $G$ is a finite cyclic group of order $n$. Then every nontrivial idempotent element in $\mathcal{K}G$ is an S-idempotent.

**Proof:** By [5], $\mathcal{K}G$ has $2^n - 2$ nontrivial idempotent elements, let $\alpha = \sum_{i=0}^{n-1} r_i g^i \in \mathcal{K}G$ be an idempotent element.

Put $\beta = \sum_{i=0}^{n-1} (-r_i)g^i \in \mathcal{K}G$. Hence

$$\beta^2 = \left(\sum_{i=0}^{n-1}(-r)_i g^i\right)^2 = \left((-1)\sum_{i=0}^{n-1} r_i g^i\right)^2 = \sum_{i=0}^{n-1} r_i g^i = \alpha$$

Now, $\alpha\beta = \sum_{i=0}^{n-1} r_i g^i \sum_{i=0}^{n-1} (-r_i)g^i = (-1)\left(\sum_{i=0}^{n-1} r_i g^i\right)^2 = \sum_{i=0}^{n-1} (-r_i)g^i = \beta$.

Therefore every nontrivial idempotent in $\mathcal{K}G$ is an S-idempotent.

Recall that $\beta$ called Smarandache Co-idempotent of $\alpha$ [1]. The following example shows that the Smarandache co-idempotent need not be unique in general.

**Example 2.1.**

Let $G$ be a cyclic group of order 3, and $\mathcal{K}$ is an algebraically closed field of characteristic 0, and let $\alpha = \sum_{i=0}^{n-1} r_i g^i \in \mathcal{K}G$. If $\alpha$ is an idempotent element, then by [5], the values of $r_0$, $r_1$ and $r_2$ are followings

| $r_0$ | 0 | $\frac{1}{3}$ | $\frac{1}{3}$ | $\frac{1}{3}$ | $\frac{1}{3}$ | $\frac{1}{3}$ | $\frac{1}{3}$ | 1 |
|---|---|---|---|---|---|---|---|---|
| $r_1$ | 0 | $\frac{1}{3}$ | $\frac{-1+\sqrt{3}\,i}{6}$ | $\frac{-1+\sqrt{3}\,i}{6}$ | $\frac{1}{3}$ | $\frac{-1+\sqrt{3}\,i}{6}$ | $\frac{-1+\sqrt{3}\,i}{6}$ | 0 |
| $r_2$ | 0 | $\frac{1}{3}$ | $\frac{-1+\sqrt{3}\,i}{6}$ | $\frac{-1+\sqrt{3}\,i}{6}$ | $\frac{1}{3}$ | $\frac{-1+\sqrt{3}\,i}{6}$ | $\frac{-1+\sqrt{3}\,i}{6}$ | 0 |

Consider the S-idempotents,
$\alpha_1 = \frac{2}{3} - \frac{1}{3}g - \frac{1}{3}g^2$, $\alpha_2 = \frac{2}{3} + \frac{1+\sqrt{3}\,i}{6}g + \frac{1-\sqrt{3}\,i}{6}g^2$ and $\alpha_3 = \frac{2}{3} + \frac{1-\sqrt{3}\,i}{6}g + \frac{1+\sqrt{3}\,i}{6}g^2$. For each ($1 \leq i \leq 3$), $\alpha_i$ has three Co-idempotents we denote them by $\beta_{ij}$ ($1 \leq j \leq 3$). They are $\beta_{11} = \frac{-2}{3} + \frac{1}{3}g + \frac{1}{3}g^2$, $\beta_{12} = \frac{\sqrt{3}\,i}{3}g - \frac{\sqrt{3}\,i}{3}g^2$, $\beta_{13} = \frac{-\sqrt{3}\,i}{3} + \frac{\sqrt{3}\,i}{3}g$, $\beta_{21} = \frac{-2}{3} - \frac{1-\sqrt{3}\,i}{6}g + \frac{-1+\sqrt{3}\,i}{6}g^2$,

$\beta_{22} = \frac{-3+\sqrt{3}\,i}{6}g + \frac{-3-\sqrt{3}\,i}{6}g^2$, $\beta_{23} = \frac{3-\sqrt{3}\,i}{6}g + \frac{1+\sqrt{3}\,i}{6}g^2$, $\beta_{31} = \frac{-2}{3} - \frac{1-\sqrt{3}\,i}{6}g - \frac{1+\sqrt{3}\,i}{6}g^2$,
$\beta_{32} = \frac{-3-\sqrt{3}\,i}{6}g + \frac{3+\sqrt{3}\,i}{6}g^2$, $\beta_{33} = \frac{3+\sqrt{3}\,i}{6}g + \frac{-3-\sqrt{3}\,i}{6}g^2$, respectively. We see that $\alpha_1\beta_{1j} = \beta_{1j}$, $\alpha_2\beta_{2j} = \beta_{2j}$ and $\alpha_3\beta_{3j} = \beta_{3j}$, $\beta_{1j}^2 = \alpha_1$, $\beta_{2j}^2 = \alpha_2$ and $\beta_{3j}^2 = \alpha_3$, for each ($1 \leq i \leq 3$).

**Theorem 2.2.**

Let $\mathcal{K}$ b an algebraically closed field of characteristic 0 and $G = \mathbb{Z}_m \times \mathbb{Z}_n$. Then every nontrivial idempotent element in $\mathcal{K}G$ is an S-idempotent.

**Proof:** If $m, n$ are relatively prime, then the proof is given in Theorem 2.1, since $\mathbb{Z}_m \times \mathbb{Z}_n \cong \mathbb{Z}_{mn}$ is cyclic. If $m$ and $n$ are not relatively prime, for each $(k,j) \in G$ let $(k,j) = g_{kn+j}$ ($0 \leq k \leq m-1$, $0 \leq j \leq n-1$, and let

$\alpha = \sum_{i=0}^{mn-1} r_i g_i \in \mathcal{K}G$ be an idempotent element [6]. Take $\beta = \sum_{i=0}^{mn-1}(-r_i) g_i \in \mathcal{K}G$, then it is clear that

$$\beta^2 = \alpha \text{ and } \alpha\beta = \beta.$$

Therefore every idempotent element in $\mathcal{K}G$ is an S-idempotent.

Finally we concern the group ring $\mathcal{R}G$ where $\mathcal{R}$ is an integral domain and $G$ is a finite group of order $n$. We give a condition under which $\mathcal{R}G$ contains S-idempotents.

**Theorem 2.3.**

Let $\mathcal{R}$ be an integral domain, and let $G$ be a finite group of order $n$. If some prime divisor $p$ of $n$ is a unit in $\mathcal{R}$ and
1) $p^3 = p^{-1}$ or
2) $p = p^{-1}$ or
3) $p = 2$.

Then the group ring $\mathcal{R}G$ has S-idempotent.

**Proof:** 1) Since $p$ is a prime dividing $n$, and $p$ is a unit in $\mathcal{R}$ then by [7] $\alpha = p^{-1}\sum_{x \in H} x$ is a nontrivial idempotent where $\mathcal{H}$ is a subgroup of $G$ of order $p$. Let $\beta = p\sum_{x \in H} x$. Then
$\alpha\beta = p^{-1}p \sum_{x \in H} x \sum_{x \in H} x = p \sum_{x \in H} x = \beta$,
and $\beta^2 = p^2(\sum_{x \in H} x)^2 = p^3 \sum_{x \in H} x = p^{-1} \sum_{x \in H} x = \alpha.$

Hence $\alpha$ is a S-idempotent.

2) we have $\alpha = p^{-1} \sum_{x \in H} x$ is a nontrivial idempotent. Let $\beta = \sum_{x \in H} x$. Then

$\alpha\beta = p^{-1} \sum_{x \in H} x \sum_{x \in H} x = \sum_{x \in H} x = \beta$,
and $\beta^2 = (\sum_{x \in H} x)^2 = p \sum_{x \in H} x = p-1 x\in H x = \alpha.$

Therefore $\alpha$ is a S-idempotent.

3) Since $p = 2$ divides $n$, then $|G| = 2k$ and $\alpha = 2^{-1}(1 + g^k)$. Let $\beta = (1 + g^k) - \alpha$. Then it is clear that $\beta^2 = \alpha$ and $\alpha\beta = \beta$. So $\alpha$ is an S-idempotent.